\documentclass[a4paper,10pt]{article}

\usepackage{alltt}
\usepackage{graphicx}
\usepackage{hyperref}

\usepackage{amssymb}
\usepackage{amsmath}
\usepackage{theorem}
\usepackage{graphics}
\usepackage{epsfig}

\newtheorem{theorem}{Theorem}[section]
{\theorembodyfont{\rmfamily}}
\newtheorem{proposition}[theorem]{Proposition}
\newtheorem{corollary}[theorem]{Corollary}
\newtheorem{conjecture}[theorem]{Conjecture}
\newtheorem{lemma}[theorem]{Lemma}
\newtheorem{remark}[theorem]{Remark}
{\theorembodyfont{\rmfamily}\newtheorem{example}[theorem]{Example}}
{\theorembodyfont{\rmfamily}}

\newenvironment{Proof}
{\begin{trivlist}\item[]{{\sc Proof.}}}{\hfill{$\square$}\noindent\end{trivlist}}

\usepackage{color}

\def\notsuccsim{\succsim\kern-11pt/\kern3pt}
\def\notsuccsimi{\succsim\kern-9pt/\kern2pt}
\def\notsuccsimii{\succsim\kern-10pt/\kern2pt}
\def\notsucc{\succ\kern-9pt/\kern3pt}
\def\notB{B\kern-12pt/\kern3pt}
\def\card#1{\vert#1\vert}

\begin{document}

\title{
\textbf{\Large Enumeration of weighted games with minimum and an analysis of voting power for bipartite complete games with minimum}}

\author{
Josep Freixas\footnote{
Department of Applied Mathematics III and High Engineering School
(Manresa Campus), Technical University of Catalonia (Spain).
e--mail: josep.freixas@upc.edu. Research partially funded by Grants
SGR 2009--1029 of
\emph{Generalitat de Catalunya} and MTM 2012--34426 from the Spanish Economy
and Competitiveness Ministry, from the Spanish Science and Innovation Ministry. E-mail: {josep.freixas@upc.edu}} \  {\normalsize and} \ Sascha Kurz\footnote{Department of Mathematics, Physics and Computer Science, University of Bayreuth, 95440 Bayreuth, Germany
Tel.: +49-921-557353   Fax: +49-921-557352.  E--mail: sascha.kurz@uni-bayreuth.de}
}
\date{\small\today}

\maketitle

\begin{abstract} This paper is a twofold contribution. First, it contributes to the problem of enumerating some classes of simple games and in particular provides the number of weighted games with minimum and the number of weighted games for the dual class as well. Second, we focus on the special case of bipartite complete games with minimum, and we compare and rank these games according to the behavior of some efficient power indices of players of type 1 (or of type 2). The main result of this second part establishes all allowable rankings of these games when the Shapley-Shubik power index is used on players of type 1.

\vskip 3mm \noindent \textit{Key words}: simple game; weighted and complete games; enumerations; Shapley-Shubik power index; Banzhaf power indices.
\vskip 3mm

\noindent {\it Math. Subj. Class. (2000)}: Primary 91A12, 91A40,
91A80, 91B12.

\noindent {\it JEL Class.}: C71, D71.
\end{abstract}

\section{Introduction}

The study of voting systems can be traced back to the late nineteenth century, when Dedekind
studied monotonic Boolean functions. In the context of voting systems
these functions correspond to simple games. In their seminal book, von Neumann and Morgenstern~\cite{vNMo44} came up with the definition of a simple game as a type of cooperative game where the payoffs to coalitions are either $1$ or $0$, i.e., coalitions can be considered either winning or losing.

A particular case of simple games, and possibly the most important subcase, is that of weighted games, in which weights are assigned to players and a threshold is set so that a coalition is winning if and only if the sum of weights of its players is at least the threshold. This is natural in Parliaments and also in corporate voting when different shareholders may own different numbers of shares. Two natural extensions of weighted games have also been thoroughly studied: (1) complete games and (2) simple games with small dimension. In this paper we deal with a particular class of complete games, the so called ``complete games with minimum" (see, e.g.,~\cite{FrPu98} and~\cite{FrPu08}).

It turns out that every complete game with five or fewer players is weighted, so the smallest possible illustrations of complete non-weighted games occur for six players: $y_1,y_2,b_1,b_2,b_3,$ and $b_4$ (here $y$ means players of yellow type, whereas $b$ means players of blue type), and we declare that a coalition is winning if and only if it contains: at least three players and at least one of them is yellow. Intuitively, it is clear that all the yellow players have the same influence (according to the desirability relation), and all the blue players have the same influence, but the yellow players have more influence than the blue players --suggesting a complete (weak) ordering for the players in this example of a voting system. In terms of the language we introduce later (in Section 2) this simple game is complete but not weighted.

Note that e.g., the coalitions of type $\{y_1,y_2,b_i\}$ for $i=1,2,3,4$ are minimal winning since all players contained are essential for the coalition to be winning. The same occurs for the coalitions $\{y_i,b_j,b_k\}$ for $i=1,2$ and $1 \leq j < k \leq 4$. However, this latter set of coalitions has an additional singularity: none of the players in these coalitions can be replaced by a weaker player. E.g., we cannot replace in these coalitions the yellow player for a blue player since the new coalition obtained would not be winning. In terms of the language we introduce later (in Section 2) we say that the coalitions of type $\{y_i,b_j,b_k\}$ for $i=1,2$ and $1 \leq j < k \leq 4$ are shift-minimal winning coalitions. On the contrary, if in a shift-minimal winning coalition we replace a weaker player by a stronger one we obtain a minimal, but not shift-minimal, winning coalition. Finally, observe that all shift-minimal winning coalitions have the same number of players of each color, i.e., they all contain one yellow player and two blue players and this information can be encapsulated in the vector: (1,2) where the first component represents the number of yellow players and the second the number of blue players. Then, we refer to the game as being complete with only one type of shift-minimal winning coalitions, or equivalently, a complete game with minimum as was denominated in~\cite{FrPu98} and~\cite{FrPu08}. Note that in the previous example there is a bipartition between types of players: yellow players and blue players. Hence, the example introduced is a bipartite complete game with minimum.

The first part of the paper deals with enumerations for weighted games with minimum, while the second part deals with rankings of players for power indices in bipartite complete games with minimum. The dimension of complete games with minimum is studied in~\cite{FrPu08}. For instance, the previous example has dimension 2 and, therefore, it decomposes as the intersection of two weighted games: $[5;3,3,1,1,1,1] \cap [3;1,1,1,1,1,1]$ (the notation for a representation of a weighted game is introduced in the preliminaries section). Most existing voting systems have a small dimension. E.g., the current voting system of the European Council is an example of a complete game with dimension 3, i.e., it decomposes as an intersection of three weighted games which cannot be simplified to an intersection of fewer weighted games~\cite{Fre04}.

The voting system to amend the Canadian Constitution is an example of a
non-weighted game which meets both requirements: it is complete (and has only one type of shift-minimal winning coalitions) and has dimension $2$, i.e., it decomposes as the intersection of two weighted games. Since 1982, an amendment to the Canadian Constitution can become law only if it is approved by at least seven of the ten Canadian provinces, subject to the proviso that the approving provinces have, among them, at least half of Canada's population. It was first studied in Kilgour~\cite{Kil83}. A census (in percentages) taken from 1960 for the Canadian provinces was: Prince Edward Island ($1\%$), Newfoundland ($3\%$), New Brunswick ($3\%$), Nova Scotia ($4\%$), Manitoba ($5\%$), Saskatchewan ($5\%$), Alberta ($7\%$), British Columbia ($9\%$), Quebec ($29\%$) and Ontario ($34\%$). This is another example of a bipartite complete game with minimum and the vector representing all shift-minimal winning coalitions is: (1,6) where the first component indicates that exactly one of the two most populated provinces votes in favor of the voted law and 6-out-of-8 of the other provinces vote in favor of the voted law as well. Games of this type are the object of study in this paper.

This paper primarily concerns enumerations. The number of complete games is known up to nine players only~\cite{FrMo10OMS}, and the number of weighted games is also known up to nine players~\cite{Kur12}. A seminal result on enumeration formulas for weighted games and complete games is May's theorem~\cite{May52}, and many other results have followed, e.g., the enumeration of weighted games with up to six players dates back at least to 1962~\cite{MTK62}.

The mathematical structure of complete games was studied in detail in the nineties by several scholars, e.g., in~\cite{KrSu95} and~\cite{CaFr96}. In the latter work a system of quantities (called characteristic invariants) is associated with every complete game and their basic properties are stated. It is shown that these quantities determine the game (uniqueness) and that every such system is associated with some complete game (existence).

According to this classification the simplest case arises when the matrix (one of the two components of the characteristic invariants) has only one shift-minimal winning vector (which corresponds therefore to a set of closely related shift-minimal winning coalitions that are enough to generate the complete game). These games have been studied in~\cite{Fre97} and~\cite{FrPu98}. The first paper provides necessary and sufficient conditions to determine whether a game of this type is weighted. In the second paper the characteristic invariants are used to ease the calculus of different types of solutions of the game like the nucleolus, the kernel and semivalues.

The interest for this type of structures has also emerged in the field of Cryptography.
The access structure in a secret sharing scheme (see e.g.,~\cite{Sti92}) can also be modeled by a simple game. To this end Simmons~\cite{Sim90} introduced the concept of a hierarchical access structure. Such an access structure stipulates that agents are partitioned into $m$ levels, and a sequence of thresholds $k_1<k_2<\dots<k_m$ is set, so that a coalition is authorized if and only if it has $k_1$ agents of the first level and $k_2$ agents of the first two levels and $k_3$ agents of the first three levels etc. These hierarchical structures are called \emph{conjunctive} since all the $m$ conditions must be satisfied for a coalition to be authorized. If only one of the $m$ conditions must be satisfied for a coalition to be authorized, then the structure is called \emph{disjunctive}. A typical example of a conjunctive hierarchical game would be the United Nations Security Council, where for the passage of a resolution all five permanent members must vote for it \emph{and} also at least nine members in total. The ideality of disjunctive games was proved by Brickell~\cite{Bri89}, while the ideality of conjunctive games was proved by Tassa~\cite{Tas07}. Ideality means they can carry the most informationally efficient secret sharing scheme and be completely secure (i.e., not giving any information about the secret to unauthorized coalitions). Gvozdeva et al.~\cite{GHS12} relate these two types of structures with complete games with one shift-minimal winning vector and with complete games with one-shift maximal losing vector.

In this paper we use game theoretic methods and terminology, and we talk about complete games with minimum or, equivalently, complete games with a unique shift-minimal winning vector (instead of hierarchical conjunctive structures) and games with a unique shift-maximal losing vector (instead of hierarchical disjunctive structures).

Here we enumerate weighted games with minimum. We find a polynomial formula as a function of the number of players of the game. This complements the corresponding known result for the enumeration of complete games of this type~\cite{FrPu98}.

As for the second contribution, we recall that the distribution of power in some important real-world institutions (the International Monetary Fund, the voting system of the World Bank, the United Nations Security Council, the procedure to amend the Canadian Constitution, etc.) has been extensively studied, e.g., in~\cite{Lee02a},~\cite{Lee02b},~\cite{AlBo05},~\cite{TaPa08},~\cite{FeMa98},~\cite{Str82} and~\cite{Kil83} to cite just some references.

We consider here the set of bipartite complete games with minimum, i.e., complete games with two types of equivalent players and one shift-minimal winning vector, and discuss the possible rankings of these games given by the Shapley-Shubik power index for a player of type 1. The main result of this part establishes all the allowable rankings for the power of players of a given type in bipartite complete games with minimum for which the number of players of each type is fixed. We do remark that many papers in the literature have been devoted to study whether two or more power indices provide the same rankings in each game (see e.g.,~\cite{DiMo02},~\cite{CaFr08},~\cite{Fre10IJGT}). However, as far as we know, very little has been done on comparing power of players in different games.

The paper is organized as follows.  In Section 2 we review some basic concepts and definitions of simple games, revise the terminology of the characteristic invariants for complete games and recall the known enumerations for complete games and for weighted games. In Section 3 we obtain a formula for the number of
weighted games with one shift-minimal winning vector and deduce some consequences. In Section 4 we do a comparison of power for different complete games with two types of equivalent players and one shift-minimal winning vector. We prove that a limited number of rankings are possible for the Shapley-Shubik power index and we formulate a similar conjecture for the relative Banzhaf index. A study of duality in Section 5 permits us to extend the results obtained in the two previous sections to complete games with one shift-maximal losing vector. Some hints for future research are given in Section 6.

\section{Preliminaries}

This preliminary section is organized into five subsections.
The first two refer to simple games in general and complete games in particular. The remaining three recall a result on the structure of complete games that will be essential for our purposes, previous results found in the literature on enumerations of games, and some power indices.

\subsection{Simple games}

A (monotonic) \emph{simple} game is a pair $(N,\mathcal{W})$ where
$N=\{1,2,...,n\}$ and $\mathcal{W}$ is a collection of
subsets of $N$ such that:
\begin{itemize}
\item[$i)$] $\emptyset \notin \mathcal W$,
\item[$ii)$] $N \in \mathcal W$,
\item[$iii)$] if $S\in \mathcal{W}$ and $S\subseteq T$, then $%
T\in \mathcal{W}$.
\end{itemize}
From now on we will omit the term monotonic. Simple games can be viewed as models of voting
systems in which a single alternative, such as a bill or an
amendment, is pitted against the status quo. The set $N$ is called
the \emph{grand coalition}, its
members are called \emph{players} and its subsets \emph{coalitions}, and the subsets in $%
\mathcal{W}$ are called \emph{winning coalitions}. The intuition here is
that a set $S$ is a winning coalition if and only if the bill or amendment
passes when the players in $S$ are precisely the ones who vote for it. A
subset of $N$ that is not in $\mathcal{W}$ is called a \emph{losing
coalition} and the collection of losing coalitions is denoted by ${\mathcal L}%
.$ If each proper subcoalition of a winning coalition is losing,
this winning coalition is called \emph{minimal}. The set of
minimal winning coalitions is denoted by $\mathcal{W}^{m}$. It
should be noted that a simple game is completely determined by its
minimal winning coalitions. If each proper coalition containing a
losing coalition is winning, this losing coalition is called
\emph{maximal}. The set of maximal losing coalitions is denoted
by ${\mathcal L}^{M}$ and it also determines the game.

Let $(N,\mathcal W)$ be a simple game. The \emph{dual} game of $(N, \mathcal W)$  is the game $(N, \mathcal W^{\ast})$ where
$\mathcal W^{\ast} = \{ S \subseteq N \, : \, N \setminus S \notin \mathcal W \}$.

A player $i$ has \emph{veto} in a simple game $(N,\mathcal{W})$ if $S\in
\mathcal{W}$ implies $i\in S$. A player $i\in N$ is called a \emph{null player} in $(N,%
\mathcal{W})$ if $i\notin S$ for every $S\in \mathcal{W}^{m}.$ A player $%
i\in N$ is a \emph{dictator} if and only if $\mathcal{W}^m = \{\{i\}\},$ in which case the
remaining players in $N$ become null players. Note that a dictator is the most extreme form
of having veto.

A simple game $(N,\mathcal{W})$ is a \emph{weighted
game} if it admits a representation by means of $n$
non-negative real numbers $w_1, \dots,w_n$ and a positive real number $q$ such that $S\in
\mathcal{W}$ if and only if $w(S)\geq q,$ where $w(S)=\sum\limits_{i\in
S}w_{i}$ for each coalition
$S\subseteq N$. The number $q$ is called the \emph{quota} of the game and $%
w_{i} $ the \emph{weight} of player $i.$ From now $[q;w_{1},...,w_{n}]$
will mean the representation of $(N,\mathcal{W})$ by means of weights
$w_1,\dots,w_n$ and quota $q$. The weighted representation (whenever it exists) is never unique. For instance, $[c\cdot q; c \cdot w_{1}, \dots, c \cdot w_{n}]$ is also a representation of $(N,\mathcal{W})$ for all $c>0$.

Two simple games $(N,\mathcal{W})$ and $(N^{\prime },\mathcal{W}^{\prime })$
are said to be \emph{isomorphic} if there exists a bijective map $%
f:N\rightarrow N^{\prime }$ such that $S\in \mathcal{W}$ if and only if $f(S)\in
\mathcal{W}^{\prime }.$

Let $(N,\mathcal{W})$ be a simple game. Set $\mathcal{W}_{i}=\{S\in \mathcal{%
W}:i\in S\}$ and let $\tau _{ij}:N\rightarrow N$ denote the
\emph{transposition} of players $i,j\in N$ (i.e., $\tau
_{ij}(i)=j,\;\tau _{ij}(j)=i\text{ and } \tau _{ij}(k)=k$ for
$k\neq i,j$). The individual {\it desirability relation},
introduced by Isbell (see~\cite{Isb56} as well as~\cite{Isb58}) and
later generalized by Maschler and Peleg~\cite{MaPe66}, is the
binary relation $\succsim$ on $N$:
\[
i\,\succsim\,j \ \mbox{ if and only if }\ \tau _{ij}(\mathcal{W}_{j})\subseteq
\mathcal{W}_{i},
\]
meaning that $i$ is \emph{at least as desirable as} $j$ as a coalition
partner. It is easy to see that $\succsim$ is a preorder (i.e., a reflexive and transitive relation),
we abbreviate $i\,\succsim\,j$, $j\,\succsim\,i$ by $i \approx j$ and say that $i$ and $j$ are equi-desirable players ($\approx$ is an equivalence relation in $N$), and we abbreviate $i\,\succsim\,j$, $j\,\notsuccsim\,i$ by $i \succ j$ and say that $i$ is strictly more desirable than $j$ as a coalition partner.

The relation $\succsim$ induces an ordering $\geq $ in the set of $\approx$-classes $N/{\approx}=\{N_{1},...,N_{t}%
\}.$ Thus, $N_{p}\geq N_{q}$ if and only if $i\,\succsim\,j$ for any $i\in N_{p}$
and any $j\in N_{q}.$

\subsection{Complete games}

The desirability is not always complete (total). Then, if any two players are comparable by $\succsim$, $(N,\mathcal{W})$
is said to be a \emph{complete game;}\footnote{%
Complete games are also known in the literature of simple games as
linear games or directed games.} in this case, the $\approx$--classes
are linearly ordered by $\geq$ . We say that a complete game has
\emph{trivial classes} if it possesses either veto or null players.
Notice that each weighted game is complete because $w_{i}\geq w_{j}$ implies
$i\,\succsim\,j.$

A coalition $S \in \mathcal W$ is \emph{shift-minimal} winning if
$(S \setminus \{i\}) \cup \{j\} \notin \mathcal W$ for all $i \in S$
and $j \notin S$ with $i \succ j$. Note that a winning coalition
can be minimal without being shift-minimal.

\begin{example} \label{e:5players}
Let $(N, \mathcal W)$
be the simple game defined by $N=\{1,2,3,4,5\},$ and
$\mathcal W^m = \{ S \subseteq N \,  : \, \card{S}=3, \, S \neq \{3,4,5\} \}$.
It is easy to check that $N$ decomposes into a bipartition of equivalent players
$N_1= \{1,2\}$ and $N_2=\{3,4,5\}$, and $i \succ j$ for all $i\in N_1$ and $j \in N_2$. Coalitions $\{1,2,3\}, \{1,2,4\}, \{1,2,5\}$ are minimal winning but not shift-minimal winning in $(N, \mathcal W)$. The remaining six winning coalitions of cardinality 3 are shift-minimal winning.
\end{example}

From now on we only deal with complete games and without loss of generality we assume $1\succsim 2 \succsim \dots \succsim n$ in the following. We can partition
the whole set $N$ of players into equivalence classes $N_1,\dots, N_t$ and say that the complete game consists of $t$~types of (weakly) ordered players.
By $n_i$ we denote the cardinality of the set $N_i$ for $1\le i\le t$. Coalitions are categorized into different types, which can be described by a vector $(m_1,\dots,m_t)$ meaning $m_i$-out-of-$n_i$ players  (from the set $N_i$) for $1\le i\le t$.

Let us consider Example~\ref{e:5players} with $n_1=2$ and $n_2=3$. Due to the assumed ordering of the players we have $N_1=\{1,2\}$ and $N_2=\{3,4,5\}$ and $i \succ j$ for all $i=1,2$ and $j=3,4,5$ so that we can write $N_1 > N_2$. With this, the vector $(1,2)$ is the type of coalitions $\{1,3,4\}$, $\{1,3,5\}$, $\{1,4,5\}$, $\{2,3,4\}$, $\{2,3,5\}$ and $\{2,4,5\}$. Since we have $1~\approx~2$ and $3~\approx~4~\approx~5$ either all these six coalitions are winning or they are all losing and we can therefore speak of a winning or a losing vector. In Example~\ref{e:5players} $(1,2)$ is a (shift-minimal) winning vector.

Let $(N,\mathcal W)$ be a simple game and $N_h$ be the classes of equally desirable players for $1\le h\le t$. We call a vector $\widetilde{m}:=(m_1,\dots,m_t)$, where $0\le m_h\le \left|N_h\right|$ for $1\le h\le t$, a
\emph{winning vector} if $S \in \mathcal W$, where $S$ is an arbitrary coalition of $N$ containing exactly $m_h$ elements of $N_h$
for $1\le h\le t$. Analogously, we call such a vector a \emph{losing vector} if $S \in \mathcal L$, where $S$ is an arbitrary
coalition of $N$ containing exactly $m_h$ elements of $N_h$ for $1\le h\le t$.

Following~\cite{CaFr96}, several concepts of ordering among vectors (i.e., types of coalitions) in ${\mathbb{N}_0}^t$
(where $\mathbb{N}_0=\mathbb{N} \cup \{0\}$) have to be considered. For two vectors $\widetilde{a}=(a_1,\dots,a_t) \in
{\mathbb{N}_0}^t$ and $\widetilde{b} =(b_1,\dots,b_t) \in {\mathbb{N}_0}^t$, representing types of coalitions in a complete
  simple game, we write (the standard componentwise order between vectors) $\widetilde{a} \geq \widetilde{b}$ if and only if we have $a_i \geq b_i$ for all $i=1,\dots, t$. We use $\widetilde{a} > \widetilde{b}$ if $\widetilde{a} \geq \widetilde{b}$ and $\widetilde{a} \neq \widetilde{b}$. We write $\widetilde{a} \succeq \widetilde{b}$ if and only if we have
  $
    \sum\limits_{i=1}^{k} a_i \geq \sum\limits_{i=1}^{k} b_i
  $
  for all $1\le k\le t$. For $\widetilde{a} \succeq  \widetilde{b}$ and $\widetilde{a} \neq \widetilde{b}$ we use
  $\widetilde{a} \succ \widetilde{b}$ as an abbreviation and say that they are comparable vectors with vector $\widetilde{b}$ being smaller
  than vector $\widetilde{a}$. If neither $\widetilde{a} \succeq \widetilde{b}$ nor $\widetilde{b} \succeq \widetilde{a}$ holds, we write $\widetilde{a} \bowtie \widetilde{b}$
  and say that vector $\widetilde{a}$ and vector $\widetilde{b}$ are \emph{incomparable}.

As $(1,2)$ is a winning vector in Example~\ref{e:5players}, so are $(1,3)$, $(2,2)$ and $(2,3)$ because of the monotonicity property required in the definition of simple game, but also is $(2,1)$ because Example~\ref{e:5players} is a complete game. From $(1,2)$ nothing can be deduced about the vectors $(1,1)$, $(0,3)$, $(0,2)$, $(1,0)$ and $(0,1)$. However, we can check that all the coalitions associated with these vectors are losing for Example~\ref{e:5players}.

A vector $\widetilde{m}=(m_1,\dots,m_t)$ in a complete game with $t$ types of equivalent players $(N, \mathcal W)$ is a \emph{minimal winning vector}
if $\widetilde{m}$ is a winning vector and every vector $\widetilde{m}'$ with $\widetilde{m} \mathbf{'} < \widetilde{m}$ is losing. Analogously,
a vector $\widetilde{m}$ is a \emph{maximal losing vector} if $\widetilde{m}$ is a losing vector and every vector $\widetilde{m}'$ with
$\widetilde{m}\mathbf{'} > \widetilde{m}$ is winning.
Of course, a vector is shift-minimal winning (resp., shift-maximal losing) if and only if any coalition represented by the vector is shift-minimal winning (resp., shift-maximal losing).

Similarly, a \emph{shift-minimal winning vector}\footnote{In~\cite{CaFr96} they are called $\delta$-minimal winning vectors.} $\widetilde{m}$ is a winning vector such that every vector $\widetilde{m}'$ with $\widetilde{m} \mathbf{'} \prec \widetilde{m}$ is losing. Analogously,
a vector $\widetilde{m}$ is a \emph{shift-maximal losing vector} if $\widetilde{m}$ is a losing vector and every vector $\widetilde{m} \mathbf{'}$ with
$\widetilde{m} \mathbf{'} \succ \widetilde{m}$ is winning.

\begin{figure}[htp]
  \begin{center}
    \setlength{\unitlength}{1cm}
    \begin{picture}(4,7.5)
      \put(1,0){$\begin{pmatrix}0&0\end{pmatrix}$}
      \put(1,1){$\begin{pmatrix}0&1\end{pmatrix}$}
      \put(0,2){$\begin{pmatrix}1&0\end{pmatrix}$}
      \put(2,2){$\begin{pmatrix}0&2\end{pmatrix}$}
      \put(1,3){$\begin{pmatrix}1&1\end{pmatrix}$}
      \put(0,4){$\begin{pmatrix}2&0\end{pmatrix}$}
      \put(2,4){$\begin{pmatrix}1&2\end{pmatrix}$}
      \put(1,5){$\begin{pmatrix}2&1\end{pmatrix}$}
      \put(2,6){$\begin{pmatrix}2&2\end{pmatrix}$}
\put(2,7){$\begin{pmatrix}2&3\end{pmatrix}$}
\put(2,4){$\begin{pmatrix}1&2\end{pmatrix}$}
\put(2,2){$\begin{pmatrix}0&2\end{pmatrix}$}
\put(3,3){$\begin{pmatrix}0&3\end{pmatrix}$}
\put(3,5){$\begin{pmatrix}1&3\end{pmatrix}$}
      \put(1.5,0.4){\vector(0,1){0.45}}
      \put(1.3,1.4){\vector(-3,2){0.65}}
      \put(1.7,1.4){\vector(3,2){0.65}}
\put(2.7,2.4){\vector(3,2){0.65}}
\put(2.7,4.4){\vector(3,2){0.65}}
\put(1.7,5.4){\vector(3,2){0.65}}
      \put(0.7,2.4){\vector(3,2){0.65}}
      \put(2.3,2.4){\vector(-3,2){0.65}}
\put(3.3,3.4){\vector(-3,2){0.65}}
\put(3.3,5.4){\vector(-3,2){0.65}}
      \put(1.3,3.4){\vector(-3,2){0.65}}
      \put(1.7,3.4){\vector(3,2){0.65}}
      \put(0.7,4.4){\vector(3,2){0.65}}
      \put(2.3,4.4){\vector(-3,2){0.65}}
      \put(2.5,6.4){\vector(0,1){0.45}}
    \end{picture}
    \caption{The Hasse diagram for the ordering $\succeq$ of vectors on $\widetilde{n} = (2,3)$.}~\label{fig_hasse_2_1}
  \end{center}
\end{figure}
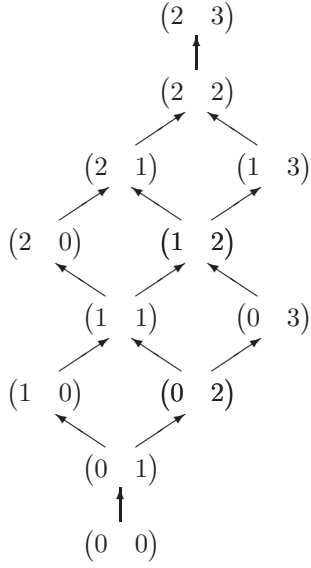

We denote by $W^m$, $W^{sm}$, $L^M$ and $L^{sM}$ the sets of minimal winning vectors, shift-minimal winning vectors, maximal losing vectors and shift-maximal losing vectors, respectively. In Example~\ref{e:5players} we have:
\begin{align}
W^m & =  \{ (2,1), \ (1,2) \}, \nonumber\\
W^{sm}& =  \{ (1,2) \}, \nonumber\\
L^M & =  \{ (2,0), \ (1,1), \ (0,3) \}, \nonumber\\
L^{sM} & =  \{ (2,0), \ (0,3) \}. \nonumber
\end{align}
The Hasse diagram for the ordering of vectors in complete games with the given hierarchy
$\widetilde{n} =(2,3)$ is shown in Figure~\ref{fig_hasse_2_1}.

\subsection{A parameterization theorem for complete games}
\label{subsec_parameterization_csg}
Carreras and Freixas have given a full parameterization of complete games, up to isomorphisms,
in~\cite{CaFr96} using vectors as models of coalitions and the partial order~$\succeq$. We denote the (decreasing) lexicographic order by $\gtrdot$, i.e.,\
we have  $(a_1,\dots,a_n) \gtrdot (b_1,\dots,b_n)$ if there is an index $1\le h\le n$ with $a_i=b_i$ for all $1\le i<h$ and $a_h>b_h$. An example is given by $(1,2,1)\gtrdot(1,1,3)$.

\vbox{
\begin{theorem}
  \label{thm_characterization_cs}

  \vspace*{0mm}

  \noindent
  \begin{itemize}
   \item[(a)] Assume that a vector $\widetilde{n} = (n_1,n_2, \dots, n_t)$ with natural coefficients and a matrix
              $$\mathcal{M}=\begin{pmatrix}m_{1,1}&m_{1,2}&\dots&m_{1,t}\\m_{2,1}&m_{2,2}&\dots&m_{2,t}\\
              \vdots&\vdots&\ddots&\vdots\\m_{r,1}&m_{r,2}&\dots&m_{r,t}\end{pmatrix}=
              \begin{pmatrix}\widetilde{m}_1\\\widetilde{m}_2\\\vdots\\\widetilde{m}_r\end{pmatrix}$$
        with natural or null coefficients are given, satisfying the following properties:
              \begin{itemize}
               \item[(i)]   $m_{1,1}>0$ and $0\le m_{i,j}\le n_j$, $m_{i,j}\in {\mathbb{N}_0}$ for $1\le i\le r$ and $1\le j\le t$,
               \item[(ii)]  $\widetilde{m}_i\bowtie\widetilde{m}_j$ for all $1\le i<j\le r$,
               \item[(iii)] for each $1\le j<t$ there is at least one row-index $i$ such that
                            $m_{i,j}>0$, $m_{i,j+1}<n_{j+1}$, and
               \item[(iv)]  $\widetilde{m}_i\gtrdot \widetilde{m}_{i+1}$ for $1\le i<r$.
              \end{itemize}
              Then, there exists a unique complete game $(N, \mathcal W)$ with invariants $\left(\widetilde{n},\mathcal{M}\right)$, i.e., with $\widetilde{n}$ as
              a vector of the cardinalities of the equivalence classes and matrix $\mathcal{M}$ where their rows consist of the shift-minimal winning vectors.
   \item[(b)] Two complete games $\left(N_1,\mathcal{W}_1\right)$ and $\left(N_2,\mathcal{W}_2\right)$
              are isomorphic if and only if $\widetilde{n}_1=\widetilde{n}_2$ and $\mathcal{M}_1=\mathcal{M}_2$.
  \end{itemize}
\end{theorem}
}

As a consequence of this theorem, any complete game can be denoted as
$(\widetilde{n}, \mathcal M)$, the pair of \emph{characteristic invariants} of the game.

\noindent
In such a vector/matrix representation (characteristic invariants) of a complete game the number of players $n$ is determined by $n=\sum\limits_{i=1}^t n_i$. Although Theorem~\ref{thm_characterization_cs} looks technical at first glance, the necessity of the required properties is easily explained. Obviously,
$n_j\ge 1$ and $0\le m_{i,j}\le n_j$ must hold for $1\le i\le r$, $1\le j\le t$. If $\widetilde{m}_i\preceq\widetilde{m}_j$ or
$\widetilde{m}_i\succeq\widetilde{m}_j$ then we have $\widetilde{m}_i=\widetilde{m}_j$ or either $\widetilde{m}_i$ or $\widetilde{m}_j$ cannot be a shift-minimal winning vector. If for a column-index $1\le j<t$ we have $m_{i,j}=0$ or $m_{i,j+1}=n_{j+1}$ for all $1\le i\le r$, then we can check whether we have $g \approx h$ for all $g\in N_j$, $h\in N_{j+1}$, which is a contradiction to the definition of the classes $N_j$ and therefore also for the numbers $n_j$. Obviously a complete game does not change if two rows of the matrix $\mathcal{M}$ are interchanged. Thus we require a given ordering of the rows to avoid repetitions: $\gtrdot$ stands for the lexicographic ordering of vectors in ${\mathbb{N}_0}^t$.

As the desirability relation is total in complete games, it defines for these games a weak ordering on the set of players.  For example, writing that a five-player complete game has hierarchy $1>2=3=4>5$ means that there is one player which has the maximum influence, another one that has the minimum influence and the other three have all the same intermediate influence, in that case we can represent the previous ordering as the vector $(1,3,1)$. We say that two complete games have the \emph{same hierarchy} if the ordering that defines the desirability relation on them is the same. Thus, if $\left(\widetilde{n}_1,\mathcal{M}_1\right)$ and $\left(\widetilde{n}_2,\mathcal{M}_2\right)$ are the characteristic invariants of two complete games, they have the same hierarchy if $\widetilde{n}_1 = \widetilde{n}_2$.

The Hasse diagram for the ordering $\succeq$ of vectors in complete games with the given hierarchy
$\widetilde{n} =(2,2)$ is shown in next Figure~\ref{fig_hasse_2_2}.
\begin{figure}[htp]
  \begin{center}
    \setlength{\unitlength}{1cm}
    \begin{picture}(3,7)
      \put(1,0){$\begin{pmatrix}0&0\end{pmatrix}$}
      \put(1,1){$\begin{pmatrix}0&1\end{pmatrix}$}
      \put(0,2){$\begin{pmatrix}1&0\end{pmatrix}$}
      \put(2,2){$\begin{pmatrix}0&2\end{pmatrix}$}
      \put(1,3){$\begin{pmatrix}1&1\end{pmatrix}$}
      \put(0,4){$\begin{pmatrix}2&0\end{pmatrix}$}
      \put(2,4){$\begin{pmatrix}1&2\end{pmatrix}$}
      \put(1,5){$\begin{pmatrix}2&1\end{pmatrix}$}
      \put(1,6){$\begin{pmatrix}2&2\end{pmatrix}$}
      \put(1.5,0.4){\vector(0,1){0.45}}
      \put(1.3,1.4){\vector(-3,2){0.65}}
      \put(1.7,1.4){\vector(3,2){0.65}}
      \put(0.7,2.4){\vector(3,2){0.65}}
      \put(2.3,2.4){\vector(-3,2){0.65}}
      \put(1.3,3.4){\vector(-3,2){0.65}}
      \put(1.7,3.4){\vector(3,2){0.65}}
      \put(0.7,4.4){\vector(3,2){0.65}}
      \put(2.3,4.4){\vector(-3,2){0.65}}
      \put(1.5,5.4){\vector(0,1){0.45}}
    \end{picture}
    \caption{The Hasse diagram for the ordering $\succeq$ of vectors on $\widetilde{n} = (2,2)$.}~\label{fig_hasse_2_2}
  \end{center}
\end{figure}
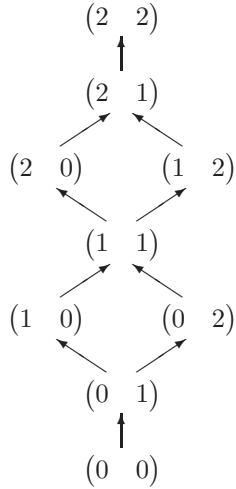

We would like to remark that for $t=1$ only $r=1$ is possible and the requirements in Theorem~\ref{thm_characterization_cs} reduce to $1\le m_{1,1}\le n_1=n$. Also for $t=2$ one can easily give a more compact formulation for the requirements in Theorem~\ref{thm_characterization_cs}. A complete description of the possible values $n_1,n_2,m_{1,1},m_{1,2}$ corresponding to a complete game with parameters $n$, $t=2$, and $r=1$ is given by
$$
\begin{array}{l}  
  1\le n_1\le n-1,\\
  n_1+n_2=n,\\
  1\le m_{1,1} \le n_1,\\
  0\le m_{1,2}\le n_2-1.
\end{array}
$$
Two important real-world examples of voting weighted games with only one shift-minimal winning vector are (see chapter 8 in~\cite{TaPa08} for more details on these two examples): the United Nations Security Council --without taking abstention into consideration-- and the procedure to amend the Canadian Constitution. These examples have $\left(\widetilde{n}, {\mathcal M}\right) = \left((5,10),(5,4)\right)$ and $\left(\widetilde{n}, {\mathcal M}\right) = \left((2,8),(1,6)\right)$ as respective characteristic invariants.

\subsection{Known enumerations for weighted games and for complete games}
\label{S:KnownEnumerations}

Let ${wg}(n,t,r)$ be the number of weighted games with $n$ players, $t$ equivalence classes $N_1,\dots,N_t$ and $r$ shift-minimal winning vectors. Let ${wg}(n,\ast, r)$\footnote{More precisely, the notation ${wg}(n,\ast, r)$ stands for $\sum\limits_{t = 1}^n{wg}(n,t,r)$.} be the number of weighted games with $n$ players and $r$ shift-minimal winning vectors (independently of the number of equivalence classes, $t$). Let ${wg}(n,t, \ast)$ be the number of weighted games with $n$ players, and $t$ equivalence classes $N_1,\dots,N_t$ (independently of the number of shift-minimal winning vectors, $r$). We identify ${wg}(n,\ast,\ast)$, i.e., the number of weighted games with $n$ players independently of the values of $r$ and $t$, with simply $wg(n)$. Analogous notations ${cg}(n,t,r)$, ${cg}(n,\ast, r)$,
${cg}(n,t, \ast)$ and ${cg}(n)$ will be used for the respective enumerations of complete games.

The first exact counting can at least be traced back to May~\cite{May52} which establishes the number of symmetric or anonymous simple games. Any such game with $n$ players admits a weighted representation $[q;\underbrace{1,1,\dots,1}_n]$  where $q \in \{1, \dots, n \}$.

Let $wg^{sym}(n)$, $cg^{sym}(n)$ and $sg^{sym}(n)$ be respectively the number of symmetric: weighted games, complete games and, simple games with $n$ players.

\begin{theorem} \label{thm_wg_n_1}
$wg(n,1,1)=cg(n,1,1)=n=wg^{sym}(n)=cg^{sym}(n)=sg^{sym}(n)$
\end{theorem}

The number of complete games with one shift-minimal winning vector was determined in~\cite{FrPu98}, and a more refined result appears in~\cite{FrPu08} (parts 1 and 2 of the next result respectively).

\begin{theorem}\label{t:cg}
\begin{enumerate}
\item $cg(n,\ast,1)=2^n-1$,
\item $$cg(n,t,1) = \left\{
                     \begin{array}{ll}
                       {n}, & \hbox{if} \ t=1 \\ \\
                       {n+1 \choose 2t-1}, & \hbox{if} \ 2 \leq t \leq \frac{n}2+1  \\ \\
                       {0}, & \hbox{otherwise}
                     \end{array}
                   \right.
$$
\end{enumerate}
\end{theorem}

Other formulas have been obtained quite recently. In~\cite{FMR12ANOR} we can find the next enumeration, where $F(n)$ are the Fibonacci numbers, which constitute a well--known sequence of integer
numbers defined by the following recurrence relation:
$F(0)=0$, $F(1)=1$, and $F(n)=F(n-1)+F(n-2)$ for all $n>1$.
\begin{theorem} \label{t:FMR12}
$
cg(n,2,\ast) = F(n+6)-(n^2 + 4n +8).
$
\end{theorem}
Let us finally remark that a formula for $cg(n,\ast,2)$ is found in~\cite{KuTa12}.
\begin{theorem} \label{t:KuTa12}
$$
cg(n,\ast,2) = \sum\limits_{t=1}^n\! cg(n,\!t,\!2)=2\cdot\left(4^n\!+\!2^n\right)\!\cdot\!
  \left(\!6\!\cdot\!\frac{8n\!-\!4n^2\!-\!3}{n(n\!-\!3)}{2n\!-\!5\choose
n\!-\!4}\!+\!\frac{2n^2\!+\!3n\!-\!2}{(n\!+\!1)(n\!-\!2)}{2n\!-\!3\choose
n\!-\!3}\!\right).
$$
\end{theorem}

In~\cite{KuTa12} it is proved that $cg(n,t,r)$ is a quasi-polynomial in $n$, if $t$ and $r$ are given, and can therefore be automatically computed (without escaping of the problem of capacity limitations).
The main purpose of Section 3 is to determine $wg(n,\ast,1)$ as well as to obtain other related finer results.

\subsection{Power indices}

Informally a \emph{power index} is a numerical measure that
estimates the a priori capacity or influence of each player in a simple game.
Of course the notion of power is complex and has been analyzed in depth by
several authors. An interesting reference is the book by Morriss~\cite{Mor02},
which analyzes power from a philosophical point of view.

Two prominent power indices are more recognized and used than others, and both known power indices, 
and both of them are based on the notion of a swing. A coalition $S$ is a \emph{swing} for
$i \in S$ if and only if $S \in \mathcal W$ but $S\setminus \{i\} \notin \mathcal W$.
Let $c_i(N, \mathcal W)$ denote the number of swings of player $i$ in game $(N, \mathcal W)$.
Then the \emph{relative Banzhaf} index~\cite{Ban65} is defined as
$$ Bz_i(N, \mathcal W) =  \dfrac{c_i(N,\mathcal W)}{\sum\limits_{j \in N} c_j(N,\mathcal W)} $$
while the \emph{absolute Banzhaf} index~\cite{Owe78} is defined as
$$ Bz_i'(N, \mathcal W) =  \dfrac{c_i(N,\mathcal W)}{2^{n-1}}. $$
The \emph{Shapley-Shubik} index~\cite{ShSh54}, which is the restriction
of the well-known Shapley value~\cite{Sha53} for cooperative games,
can be expressed as a function of the swings as follows.
Let $s$ be the cardinality of the swing $S$ for $i$ and $c_i^s(N, \mathcal W)$
be the number of swings for $i$ for coalitions $S$ of cardinality $s$. Then,
$$ SS_i(N, \mathcal W) = \sum_{s=1}^n \dfrac{(s-1)!(n-s)!}{n!} c_i^s(N, \mathcal W). $$
This less usual formulation of the Shapley-Shubik index will be helpful
here for our purposes.
The Shapley-Shubik index of a player can be viewed as his/her expected part of a
fixed total prize, i.e., the power of a player is meant to
be the player's expected payoff. Dubey and Shapley~\cite{DuSh79} proved the
following result.
\begin{proposition} \label{p:ds79}
If $(N, \mathcal W)$ is any simple game then, for all $i \in N$, we have
\begin{enumerate}
\item[(a)] $SS_i(N, \mathcal W) = SS_i(N, \mathcal W^{\ast})$.
\item[(b)] $Bz_i(N, \mathcal W) = Bz_i(N, \mathcal W^{\ast})$ and  $Bz'_i(N, \mathcal W) = Bz'_i(N, \mathcal W^{\ast})$.
\end{enumerate}
The fact that $c_i^s(N,\mathcal W) = c_i^s(N,\mathcal W^{\ast})$ for all $s=1,2,\dots,n$ justifies (a) and implies that $c_i(N, \mathcal W)=c_i(N,\mathcal W^{\ast})$, which justifies (b), a property discovered by Dubey and Shapley~\cite{DuSh79}.
\end{proposition}

\section{Counting weighted games with minimum}

If we consider a subgame of a weighted game when this game is stripped of its null and veto players, then the original game is weighted if and only if the subgame is weighted and, similarly, the original game is complete if and only if the subgame is complete, see e.g.,~\cite{TaZw99}.

Let $\left(\widetilde{n},\mathcal{M}\right)$ be a complete game with $t$ equivalence classes. By $\left(\widetilde{n},\mathcal{M}\right)\!\!\downarrow$ we denote
$$
    \left(\begin{pmatrix} n_s,&\dots,&n_e \end{pmatrix},\begin{pmatrix} m_{1,s}&\dots&m_{1,e}\\ \vdots&\ddots&\vdots\\
    m_{r,s} &\dots&m_{r,e} \end{pmatrix}\right),
$$
where $s=1$ if there are no veto players, $s=2$ if there are veto players (which would form the strongest class), $e=t$ if there are no null players, and $e=t-1$ if there are null players (which would form the weakest class). If all players of $\left(\widetilde{n},\mathcal{M}\right)$ have veto or are null players $\left(\widetilde{n},\mathcal{M}\right)\!\!\downarrow$ is  \emph{empty}.

For instance, $\left(\widetilde{n},\mathcal{M}\right)\!\!\downarrow$ does not change for the system to amend the Canadian Constitution, i.e., $\left(\widetilde{n},\mathcal{M}\right)\!\!\downarrow =
\left(\widetilde{n},\mathcal{M}\right) = \left( (2,8), (1 \, \, 6) \right)$, whereas $\left(\widetilde{n},\mathcal{M}\right)\!\!\downarrow$ reduces to
$\left((10),(4)\right)$ for the United Nations Security Council, because the five permanent members have veto right.

\begin{lemma}
  \label{lemma_count_trivial}
  Let $\widehat{wg}(n,t,r)$ be the number of non-trivial weighted games with $t$ equivalence classes $N_1,\dots,N_t$ and $r$ shift-minimal
  winning vectors. For $r>1$ or $t>2$ we have
  \begin{equation}
    \label{eq_non_trivial}
     wg(n,t,r)=\widehat{wg}(n,t,r)+\sum_{h=1}^{n-1} 2\cdot \widehat{wg}(n-h,t-1,r)+(h-1)\cdot \widehat{wg}(n-h,t-2,r),
  \end{equation}
where we define $\widehat{wg}(n,t,r)=0$ for the non-feasible cases $n<t$ or $t<1$. For $r=1$ and $t=1$ an additional $1$ has to be added to the right hand side of Equation~(\ref{eq_non_trivial}), and for $r=1$ and $t=2$ an additional term $n-1$ has to be added to the right hand side of Equation~(\ref{eq_non_trivial}).
\end{lemma}
\begin{Proof}
Every weighted game arises from a non-trivial weighted game or an empty game by appending $h_1\ge 0$
  veto players and $h_2\ge 0$ null players.
\end{Proof}

For $r=1$ and arbitrary $t$ the set of maximal losing vectors of complete games without null players, i.e., with $m_{1,t}\ge 1$, was analytically given in~\cite{FrPu08} (see next lemma).
As for $r=1$ the first indices in vector $\widetilde{m}_1$ do not carry any information we omit them.
\begin{lemma}
  \label{lemma_maximal_losing_vectors_r_1}
  For a complete game
  $$
  \left(\begin{pmatrix}n_1,\dots,n_t\end{pmatrix},\begin{pmatrix}m_{1} \dots m_{t}\end{pmatrix}\right)
  $$
  without null players, i.e., with $m_{t}\ge 1$, the complete set of shift-maximal losing vectors is given in the following matrix:
  $$
    \begin{pmatrix}
      a_{1,1} & n_2 & n_3 & \dots & n_t\\
      a_{2,1} & a_{2,2} & n_3 & \dots & n_t\\
      a_{3,1} & a_{3,2} & a_{3,3} & \dots & n_t \\
      \vdots & \vdots & \vdots & \ddots & \vdots \\
      a_{t-1,1} & a_{t-1,2} & a_{t-1,3} & \dots & n_t\\
      a_{t,1} & a_{t,2} & a_{t,3} & \dots & a_{t,t}
    \end{pmatrix}=:\begin{pmatrix}\widetilde{a}_1\\\vdots\\\widetilde{a}_t\end{pmatrix},
  $$
  where
  $$
    a_{i,j}:=\max\left\{0,\min\left\{n_j,\,-1+\sum_{h=1}^i m_{h}-\sum_{h=1}^{j-1}n_h\right \}\right \}
  $$
  for $1\le j\le i\le t$.
\end{lemma}

Note that each vector $\widetilde{a}_i$ for $1\leq i \leq t$ represents shift-maximal losing coalitions.
Indeed, these coalitions contain $\sum_{j=1}^i m_{j}-1$ strongest players (according to the desirability relation) and additionally they also contain the $\sum_{j=i+1}^t n_{j}$ weakest players, i.e.,  those players that belong to the $t-i$ weakest equivalence classes, from the $i+1$th class to the $t$th last class. If in one of these coalitions an additional player was added or a weaker player was replaced by a stronger one, then the new vector $\widetilde{b}$ representing the coalition would contain at least
$\sum_{j=1}^i m_{j}$ players for each $1\leq i \leq t$ and therefore would be winning.

The conditions~(i)--(iv) of Theorem \ref{thm_characterization_cs}-(a) reduce to $1\le m_{1}\le n_1$, $0\le m_{t}\le n_t-1$, and $1\le m_{h}\le n_h-1$ for all $h$ such that $2\le h\le t-1$. We remark that if a complete game with $r=1$ has null players, then the complete set of shift-maximal losing vectors is given by the vectors in Lemma \ref{lemma_maximal_losing_vectors_r_1} except for the last vector $\widetilde{a}_t$. We would like to remark too that $\widetilde{a}_t$ is a losing
vector which is not maximal and that there is a typo in the remark of \cite{FrPu08}, i.e., the first row (instead of the last one) must be deleted.

Having an analytic description of the shift-maximal losing vectors at hand it is not too hard to characterize the set of weighted games analytically. Indeed, already in \cite{Fre97} the non-trivial weighted games with $r=1$, i.e., those having exactly one shift-minimal winning vector, were completely classified.

\begin{theorem}
  A non-trivial complete game $\left(\widetilde{n},\mathcal{M}\right)$ with $r=1$ is a weighted game if and only if 
  either $t=1$ or $t=2$ and $m_{2}\in\{1,n_2-1\}$.
\end{theorem}
\begin{Proof}
  Due to Theorem~\ref{thm_characterization_cs} we can assume $t\ge 2$. For $t=2$ the shift-maximal losing vectors are given by
  $$ (m_{1}-1, n_2) \quad\text{and}\quad (m_{1}+c, m_{2}-1-c),$$
  where $c=\min\left(n_1-m_{1},m_{2}-1\right)$. Choosing $w_1=1$ we conclude
  \begin{eqnarray*}
    m_{1}+w_2m_{2}&>&m_{1}-1+w_2n_2,\\
    m_{1}+w_2m_{2}&>& m_{1}+c +w_2\left(m_{2}-1-c\right),
  \end{eqnarray*}
  which is equivalent to $w_2<\frac{1}{n_2-m_{2}}$ and $(1+c)w_2>c$. If $m_{2}=n_2-1$ or $c=0$, which is equivalent to $m_{2}=1$,
  there exists a solution for $w_2$. If $m_{2}\le n_2-2$ and $c\ge 1$ the two inequalities are contradicting.

  Now we consider the remaining cases $t\ge 3$. Here the vectors
  $$(m_{1}-1,n_2,\dots ,n_t ), \quad \text{ and } \quad
  (m_{1}+1,m_{2}-1,m_{3}-1,n_4,\dots,n_t)$$
  are losing vectors (the second not necessarily shift-maximal). Thus for $w_1=1$ we have
  \begin{eqnarray*}
    m_{1}+\sum_{j=2}^{t} w_j m_{j} &>& m_{1}-1+\sum_{j=2}^{t} w_j n_j,\\
    m_{1}+\sum_{j=2}^{t} w_j m_{j} &>& m_{1}+1+w_2(m_{2}-1)+w_3(m_{3}-1)+\sum_{j=4}^{t} w_j n_j,
  \end{eqnarray*}
  from which we conclude $1>w_2+w_3$ and $w_2+w_3>1$, which is a contradiction.
\end{Proof}

Thus together with Lemma~\ref{lemma_count_trivial} and Theorem~\ref{thm_wg_n_1} we can conclude:
\begin{theorem} \label{t:weightedformulas} For $n\ge 1$ we have
  \begin{eqnarray*}
    \widetilde{wg}(n,1,1) &=& n-1,\\
    \widetilde{wg}(n,2,1) &=& \left\{\begin{array}{ll}0&\text{if }n\le 2,\\n^2-6n+9&\text{if }n\ge 3,\end{array}\right.\\
    \widetilde{wg}(n,t,1) &=& 0\quad\text{if }t\ge 3,\\
    wg(n,1,1) &=& n,\\
    wg(n,2,1) &=& \left\{\begin{array}{ll}n-1&\text{if }n\le 2,\\2(n-2)^2+2&\text{if }n\ge 3,\end{array}\right.\\
    wg(n,3,1) &=& \left\{\begin{array}{ll}0&\text{if }n\le 3,\\\frac{5n^3-48n^2+157n-174}{6}&\text{if }n\ge 4,\end{array}\right.\\
    wg(n,4,1) &=& \left\{\begin{array}{ll}0&\text{if }n\le 5,\\\frac{n^4-16n^3+95n^2-248n+240}{12}&\text{if }n\ge 6,\end{array}\right.\\
    wg(n,t,1) &=& 0\quad\text{if }t\ge 5.\\
  \end{eqnarray*}
\end{theorem}

Adding up the previous enumerations we get the following compact expression for the total number of weighted games with a unique shift-minimal winning vector.

\begin{corollary} \label{c:compactwg(n,1)}
$$ {wg}(n,\ast,1) = \left\{
                        \begin{array}{ll}
                          2^n-1 & \hbox{if} \; n \leq 5 \\
                          \dfrac{n^4-6n^3+23n^2-18n+12}{12} & \hbox{if} \; n \geq 6 \\
                        \end{array}
                      \right.
    $$
\end{corollary}
Thus, the numbers of weighted games ${wg}(n,\ast,1)$ and complete games ${cg}(n,\ast,1)$ (see Theorem~\ref{t:cg}) coincide for $n\le 5$ players, but their ratio converges to zero as $n$ increases. An asymptotic upper bound for weighted games is given in~\cite{KKZ10} and an asymptotic lower bound for complete games is given in~\cite{PeSi85}, where these games are called regular Boolean functions. Many useful and accurate asymptotic estimations for simple games and subclasses of them are provided in~\cite{Kor03}.

\section{Analysis of voting power for bipartite complete games with minimum}

Many papers have been devoted to study classes of games for which two or more power indices provide the same rankings in every game of such classes, see e.g.,~\cite{DiMo02},~\cite{CaFr08},~\cite{Fre10IJGT}. However, as far as we know, very little has been studied about comparisons of different games depending on how a given power index acts on them.

In this section we will consider complete games with two types of equivalent players and with only one shift-minimal winning vector. For these games we will give explicit formulas to calculate some efficient power indices. Since the games have only two types of equivalent players and the considered power indices are efficient, it will be possible to totally rank these games  according to the behavior of the power index on players that belong to a fixed class. Of course, the order obtained considering the power over a type of players for bipartite complete games with games with minimum is reversed if the power index is evaluated over players that belong to the other equivalence class.

An important tool to this purpose is monotonicity in Young's sense~\cite{You85} which was used to give a characterization of the Shapley value that avoids additivity. If we assume a fixed set of players $N$ we write $ \mathcal W \, B_i \, \mathcal W'$ whenever if $S$ is a swing for $i$ in $\mathcal W'$ then $S$ is a swing for $i$ in $\mathcal W$. Relation $B_i$ allows us to  qualitatively compare the position of a given player $i$ in two games. However, we do not see a direct application of such a monotonicity for most of the cases we are going to study.

\subsection{Allowable rankings for the Shapley-Shubik index}

The main purpose of this subsection is to study the allowable hierarchies that the Shapley-Shubik index produces when applied to different bipartite complete ($t=2$) games with minimum ($r=1$). We also illustrate the difficulty to extend similar results for other power indices.

\begin{proposition} \label{p:SSpower}
  Let $\left(\widetilde{n},\mathcal{M}\right)$ be a complete game with $\widetilde{n}=(n_1,n_2)$ and
  $\mathcal{M}=\begin{pmatrix}a&b\end{pmatrix}$, where $a\ge 1$ and $b\le n_2-1$.
  \begin{itemize}
    \item[(1)] For a player of type~$1$ the number of coalitions where he is a swing player is given by
               $$
                 c_1=\sum_{i=b+1}^{n_2}{n_1-1\choose a-1}\cdot{n_2\choose i}\,+\,
                 \sum_{i=0}^{\min\{b,n_1-a\}}{n_1-1 \choose a+i-1}\cdot{n_2\choose b-i}.
               $$
    \item[(2)] For a player of type~$2$ the number of coalitions where he is a swing player is given by
               $$
                 c_2=\sum_{i=0}^{\min\{b-1,n_1-a\}}{n_1 \choose a+i}\cdot{n_2-1\choose b-i-1}.
                 \footnote{Unfortunately we cannot apply Vandermonde's identity $\sum\limits_{j=0}^{k}{m\choose j}{n\choose k-j}={m+n\choose k}$
                 directly.}
               $$
    \item[(3)] The Shapley-Shubik power index $SS_1(a,b,n_1,n_2)$ of a player of type~$1$ is given by
    $$
    \begin{array}{l}
      \dfrac{1}{n!}\cdot\sum\limits_{i=b+1}^{n_2}{n_1-1\choose a-1}\cdot{n_2\choose b}\cdot(a+i-1)!\cdot(n-a-i)!\,+\,   \\ \\
      \dfrac{(a+b-1)!\cdot(n-a-b)!}{n!}\cdot \sum\limits_{i=0}^{\min\{b,n_1-a\}}{n_1-1 \choose a+i-1}\cdot{n_2\choose b-i}.
    \end{array}
    $$
   \item[(4)] The Shapley-Shubik power index $SS_2(a,b,n_1,n_2)$ of a player of type~$2$ is given by
               $$
                 \dfrac{(a+b-1)!\cdot(n-a-b)!}{n!}\cdot \sum_{i=0}^{\min\{b-1,n_1-a\}}{n_1 \choose a+i}\cdot{n_2-1\choose b-i-1}.
               $$
   \item[(5)] For $b\ge 1$ we have $SS_2(a,b,n_1,n_2)-SS_2(a,b-1,n_1,n_2)=$
              $$
                \dfrac{(a+b-2)!\cdot(n-a-b)!\cdot n_1\cdot{n_1-1\choose a-1}\cdot{n_2-1\choose b-1} }{n!}>0.
              $$
   \item[(6)] For $a\ge 2$ we have $SS_2(a-1,b,n_1,n_2)-SS_2(a,b,n_1,n_2)=$
              $$
              \left\{
                \begin{array}{ll}
                  0, & \hbox{if} \; b=0 \\
                  \dfrac{(a+b-2)!\cdot(n-a-b)!\cdot (n_2-1)\cdot{n_1\choose a-1}\cdot{n_2-2\choose b-1} }{n!}, & \hbox{if} \; b \geq 1 \\
                \end{array}
              \right.
              $$
  \end{itemize}
\end{proposition}
\begin{Proof}
  \begin{itemize}
   \item[(1)] The vectors representing coalitions where a player of type~$1$ is a swing player are given by $(a,b+1)$, $(a,b+2)$, $\dots$, $(a,n_2)$
              and $(a,b)$, $(a+1,b-1)$, $(a+2,b-2)$, $\dots$, $(c,d)$ where:
    $$(c,d) = \left\{
                \begin{array}{ll}
                  (a+b,0), & \hbox{if} \; a+b \leq n_1 \\
                  (n_1,a+b-n_1), & \hbox{otherwise.}
                \end{array}
              \right.
    $$
The number of swings for a player of type $1$ for an arbitrary vector $(x,y)$ is:
${n_1-1\choose x-1} \cdot {n_2\choose y}$.
   \item[(2)] The vectors representing coalitions where a player of type~$2$ is a swing player are given by $(a,b)$, $(a+1,b-1)$, $(a+2,b-2)$, $\dots$, $(e,f)$ where:
    $$(e,f) = \left\{
                \begin{array}{ll}
                  (a+b-1,1), & \hbox{if} \; a+b \leq n_1 \\
                  (n_1,a+b-n_1), & \hbox{otherwise.}
                \end{array}
              \right.
    $$
The number of swings for a player of type $2$ for an arbitrary vector $(x,y)$ with $y>0$ is:
${n_1\choose x} \cdot {n_2-1\choose y-1}$.
\item[(3)-(4)] These results follow from the definition given in this paper for the Shapley-Shubik index and parts (1)-(2) respectively.
   \item[(5)] For $n_1-a\le b-2$ we have
   $$
   \begin{array}{l}
     \hskip -1truecm SS_2(a,b,n_1,n_2)-SS_2(a,b-1,n_1,n_2)=  \frac{(a+b-2)!(n-a-b)!}{n!} \cdot \\ \\
     \hskip -1truecm \left[\sum\limits_{i=0}^{n_1-a}{n_1\choose a+i}{n_2-1\choose b-i-1}\left(a+b-1-\frac{(n-a-b+1)(b-i-1)}{n_2-b+i+1}\right)\right]
   \end{array}
   $$
   and for $n_1-a\ge b-1$ we have:
   $$
   \begin{array}{l}
     SS_2(a,b,n_1,n_2)-SS_2(a,b-1,n_1,n_2)=\frac{(a+b-2)!(n-a-b)!}{n!}\cdot \\ \\
     \left[(a+b-1)\sum\limits_{i=0}^{b-1}{n_1\choose a+i}{n_2-1\choose b-i-1}\,-\,(n-a-b+1)\sum\limits_{i=0}^{b-2}{n_1\choose a+i}{n_2-1\choose b-i-2}\right],
   \end{array}
   $$
both of which can be simplified to the stated expression.
    \item[(6)] Similar to (5).
  \end{itemize}
\end{Proof}

\begin{corollary}
  \label{cor_pi_orderSS}
  \begin{enumerate}
  \item For a given vector $(n_1,n_2)$ let $((m_1,m_2))$ and $((m_1',m_2'))$ be two different complete simple games, i.e., $0<m_1,m_1'\le n_1$ and $0\le m_2,m_2'<n_2$.
  If $m_1\ge m_1'$ and $m_2\le m_2'$ then $SS_2(m_1,m_2) \le SS_2(m_1',m_2')$ and $SS_1(m_1,m_2) \ge SS_1 (m_1',m_2')$, where equality holds if and only if $m_2=m_2'=0$.
  \item For a given vector $(n_1,n_2)$:
  \begin{equation} \label{E:SSbounds}
   \dfrac 1n <SS_1(1,n_2-1) \leq SS_1(m_1,m_2) \leq SS_1(n_1,1) < SS_1(c,0) = \dfrac 1{n_1}
  \end{equation}
  where $c$ is any integer number between $1$ and $n_1$.

  Inequalities~\eqref{E:SSbounds} imply
  $$ \dfrac 1n >SS_2(1,n_2-1) \geq SS_2(m_1,m_2) \geq SS_2(n_1,1) > SS_2(c,0) = 0 $$
  where $c$ is any integer number between $1$ and $n_1$.
  \end{enumerate}
\end{corollary}

\bigskip

\begin{remark}
   We have implemented a computer program which can determine the Banzhaf and the Shapley-Shubik power index for bipartite complete games with minimum. For the case $n_1=3$, $n_2=7$ we have the following ordering with respect to
$SS_1$:
  \begin{eqnarray*}
    && {\bf (3,0)=(2,0)=(1,0)>(3,1)}>(3,2)>(2,1)>(3,3)> \\
    &&(2,2)>(3,4)>(1,1)>(2,3)>(3,5)>(1,2)>(2,4)> \\
    &&(3,6)>(1,3)>(2,5)>(1,4)>(2,6)>(1,5){\bf >(1,6)}
  \end{eqnarray*}
\end{remark}

Note that, in general, for an arbitrary pair $(n_1,n_2)$ the rankings of some games with respect to
$SS_1$ or $Bz_1$ (printed in bold in the previous example) are fixed. We refer to the games of type $(c,0)$ for $c>0$ which are always tied among them and situated on the top of the ranking and, oppositely, the game $(1,n_2-1)$ which is always situated at the bottom of the ranking. These extreme games are highlighted in black in the previous example for $n_1=3$, $n_2=7$.
Additionally,
Corollary~\ref{cor_pi_orderSS} provides some constraints on the rankings of the different games with respect to $SS_1$ or $Bz_1$: \begin{enumerate} \item leaving the first component fixed:
\newline
$(3,1)>(3,2)>(3,3)>(3,4)>(3,5)>(3,6);$
\newline $(2,1)>(2,2)>(2,3)>(2,4)>(2,5)>(2,6);$ and
\newline $(1,1)>(1,2)>(1,3)>(1,4)>(1,5)>(1,6).$
\item leaving the second component fixed: \newline
$(3,1)>(2,1)>(1,1)$; $(3,2)>(2,2)>(1,2)$; $(3,3)>(2,3)>(1,3)$; \newline
$(3,4)>(2,4)>(1,4)$; $(3,5)>(2,5)>(1,5)$ and $(3,6)>(2,6)>(1,6)$.
\end{enumerate}
Taking into account all these restrictions we have that for the given hierarchy $\overline n = (3,7)$ we have $13$ weighted games, three of which are of type $(c,0)$ and give the maximum value $1/n_1=1/3$ for the $SS_1$ to three most powerful players; the ranking of all other games with respect to
$SS_1$ is strict so that there are, in principle, $10!$ potential strict orderings for the power of $SS_1$ over the set of these games, but Corollary~\ref{cor_pi_orderSS} guarantees that at most $12$ of these rankings are possible.

\subsection{Comparisons for the two Banzhaf indices}

As we shall see in this section, we find for the relative Banzhaf index a similar result to that one obtained for the Shapley-Shubik index (for less than 100 players), while it fails for the absolute Banzhaf index.

The relative Banzhaf power index $Bz_1(a,b,n_1,n_2)$ of a player of type~$1$ is given by
               $$
                 \dfrac{c_1}{n_1\cdot c_1+n_2\cdot c_2}.
               $$
The relative Banzhaf power index $Bz_2(a,b,n_1,n_2)$ of a player of type~$2$ is given by
               $$
                 \dfrac{c_2}{n_1\cdot c_1+n_2\cdot c_2}.
               $$
The absolute Banzhaf power index
$Bz_1'(a,b,n_1,n_2)$ of a player of type~$1$ is given by
               $$
                 \dfrac{c_1}{2^{n-1}}.
               $$
The absolute Banzhaf power index $Bz_2'(a,b,n_1,n_2)$ of a player of type~$2$ is given by
               $$
                 \dfrac{c_2}{2^{n-1}}.
               $$
\begin{remark}
  For the Banzhaf absolute power index, Corollary~\ref{cor_pi_orderSS} is wrong, and an example is given by $n=4$, $n_1=n_2=2$
  and the games  $((2,1))$, $((1,1))$ with Banzhaf values $\left(\frac{3}{8},\frac{1}{8}\right)$, $\left(\frac{4}{8},\frac{2}{8}\right)$, respectively. Another
  example is given by $n=7$, $n_1=3$, $n_2=4$ and the games $((3,1))$ and $((2,2))$ with absolute Banzhaf indices $\left(\frac{15}{64},\frac{1}{64}\right)$,
  $\left(\frac{26}{64},\frac{10}{64}\right)$, respectively.
\end{remark}

\begin{example}[Ties]
  Let us consider the equality case in Corollary~\ref{cor_pi_orderSS} for the relative Banzhaf power index, i.e., where the powers sum up to $1$. For $m_2=m_2'=0$
  equality holds. Other examples are given by
  \begin{itemize}
   \item $n=7$, $n_1=3$, and the games $((3,3))$, $((1,2))$ with Banzhaf numerators $(5,3)$, $(20,12)$
   \item $n=8$, $n_1=2$, and the games $((2,5))$, $((1,4))$ with Banzhaf numerators $(7,5)$, $(42,30)$
   \item $n=9$, $n_1=6$, and the games $((6,2))$, $((1,1))$ with Banzhaf numerators $(4,2)$, $(12,6)$
   \item $n=13$, $n_1=5$, and the games $((5,7))$, $((2,5))$ with Banzhaf numerators $(9,7)$, $(1044,812)$
   \item $n=13$, $n_1=9$, and the games $((6,2))$, $((1,1))$ with Banzhaf numerators $(736,288)$, $(23,9)$
  \end{itemize}
  These are all examples where $(m_1,m_2)>(m_1',m_2')$ or $(m_1,m_2)<(m_1',m_2')$.  If we assume $m_1\ge m_1'$ and $m_2\le m_2'$ then there is no further example
  for $n\le 32$.
\end{example}

\begin{remark}
  For the relative Banzhaf power index Corollary~\ref{cor_pi_orderSS} is true for all $n\le 100$.
\end{remark}
This suggests to ask whether Corollary~\ref{cor_pi_orderSS} is true for the relative Banzhaf index for all $n$.

\begin{conjecture} \label{c:BzR1}
\begin{itemize}
\item[(1)] For $b\ge 1$ we have $Bz_2(a,b,n_1,n_2)-Bz_2(a,b-1,n_1,n_2)>0.$

\item[(2)] For $a\ge 2$ we have $Bz_2(a-1,b,n_1,n_2)-Bz_2(a,b,n_1,n_2) \geq 0$
              which equals zero for $b=0$ and otherwise it is positive.
\end{itemize}
\end{conjecture}
It would imply a corollary analogous to Corollary~\ref{cor_pi_orderSS}.

\begin{corollary}
  \label{cor_pi_order Bz} \label{c:BzR2}
  \begin{enumerate}
  \item For a given vector $(n_1,n_2)$ let $((m_1,m_2))$ and $((m_1',m_2'))$ be two different complete simple games, i.e., $0<m_1,m_1'\le n_1$ and $0\le m_2,m_2'<n_2$.
  If $m_1\ge m_1'$ and $m_2\le m_2'$ then $Bz_2(m_1,m_2) \le Bz_2(m_1',m_2')$ and $Bz_1(m_1,m_2) \ge Bz_1 (m_1',m_2')$, where equality holds if and only if $m_2=m_2'=0$.
  \item For a given vector $(n_1,n_2)$:
  \begin{equation}~\label{E:Bzbounds}
   \dfrac 1n < Bz_1(1,n_2-1) \leq Bz_1(m_1,m_2) \leq Bz_1(n_1,1) < Bz_1(c,0) = \dfrac 1{n_1}
  \end{equation}
  where $c$ is any integer number between $1$ and $n_1$.

  Inequalities~\eqref{E:Bzbounds} imply
  $$ \dfrac 1n >Bz_2(1,n_2-1) \geq Bz_2(m_1,m_2) \geq Bz_2(n_1,1) > Bz_2(c,0) = 0 $$
  where $c$ is any integer number between $1$ and $n_1$.
  \end{enumerate}
\end{corollary}

However, the ranking over different games for Shapley-Shubik and relative Banzhaf index are not necessarily the same, as the following example illustrates.

    \begin{example} For the case $n_1=3$, $n_2=7$ we have checked the following rankings with respect to $Bz_1$:
    \begin{eqnarray*} &&(3,0)=(2,0)=(1,0)>(3,1)>(2,1)>(1,1)>(3,2)> \\ &&(2,2)>(3,3)>(1,2)>(2,3)>(3,4)>(1,3)>(2,4)> \\ &&(3,5)>(2,5)>(1,4)>(3,6)>(2,6)>(1,5)>(1,6) \end{eqnarray*}
    This ranking of different games with vector $\widetilde{n} = (3,7)$ does not coincide with the ranking obtained for the $SS_1$ but, restricted to weighted games, it is one of the $12$ expected rankings for the $SS_1$ according to Corollary~\ref{cor_pi_orderSS}.
    \end{example}

\begin{remark}
  One might conjecture that Corollary~\ref{cor_pi_orderSS} holds for
  some effective generalized power indices as for example semivalues. 
  Semivalues for simple games (or power semiindices) are uniquely determined as those power indices that satisfy: symmetry, positivity, dummy player property and transfer (see~\cite{CFP03}) and with specific coefficients $\{p_j\}_{j=1}^{n}$ such that $\sum_{j=1}^{n} p_j{n-1 \choose j-1}=1$ and $p_j \geq 0$ for all $j$. The absolute Banzhaf power index is given by $p_j= 1/2^{n-1}$ for all $1 \le j \leq n$ and the Shapley-Shubik power index is given by $p_j = \dfrac 1{n {n-1 \choose j-1}}$.
  The (unnormalized) power, $SV$, of player $i$
  is given by
$$
    SV(i):=\sum_{i=1}^n p_{|S|} \cdot c^s_i.
$$
where $|S|=s$. 

  Let us consider the following example: $n=10$, $n_1=7$, and the games $((7,1))$ and $((1,2))$ with power index numerators $(3p_7+3p_8+p_9,p_7)$,
  $(36p_2+p_3,35p_2)$. This example is a counterexample to Corollary~\ref{cor_pi_order Bz} if and only if
  $$
     207p_2p_7+315p_2p_8+105p_2p_9-3p_3p_7<0.
  $$
  If we assume $p_j=p_{n-j-1}$ then this is equivalent to
  $$
    p_2\cdot\left(105p_0+325p_1+207p_2-3p_3\right)<0,
  $$
  and thus it is possible.
\end{remark}

\section{Results preserved by duality}

This section establishes a simple but significant result on enumerations and highlights
that the results obtained in the previous section about power indices are preserved by duality.

Let us recall that the \emph{dual} game of $(N, \mathcal W)$  is $(N, \mathcal W^{\ast})$ where:
$\mathcal W^{\ast} = \{ S \subseteq N \, : \, N \setminus S \notin \mathcal W \}$.
Then it is not difficult to check the following:
\begin{lemma} \label{l:duals}
\begin{itemize}
\item[$(i)$]
 $(N, \mathcal W)$ is weighted if and only if $(N, \mathcal W^{\ast})$ is weighted,
  and if $[q;w_1, \dots,w_n]$ is an integer representation for $(N,\mathcal W)$  then
$[T-q+1;w_1, \dots,w_n]$ is an integer representation for $(N,\mathcal W^{\ast})$ where $ T = \sum\limits_{i=1}^n w_i$ and vice versa.
\item[$(ii)$]  $i \succsim j$ if and only if $i \succsim^{\ast} j$ where $\succsim^{\ast}$ stands for the desirability relation for game $(N, \mathcal W^{\ast})$. Thus, $(N, \mathcal W)$ is complete if and only if $(N, \mathcal W^{\ast})$ is complete, and the $\approx$-classes $N_i$ for $(N, \mathcal W)$ and its ordering $N_1 >N_2 \dots > N_t$ are preserved by $\succsim^{\ast}$, i.e., $(N, \mathcal W)$ and $(N, {\mathcal W}^{\ast})$ have the same ranking $\widetilde{n}$.
\item[$(iii)$] The complete game $(N, \mathcal W)$ has the vector $\widetilde{n}=(n_1,\dots,n_t)\in\mathbb{N}_{>0}^t$ and the matrix
              $$\mathcal{M}=\begin{pmatrix}m_{1,1}&m_{1,2}&\dots&m_{1,t}\\m_{2,1}&m_{2,2}&\dots&m_{2,t}\\
              \vdots&\vdots&\ddots&\vdots\\m_{r,1}&m_{r,2}&\dots&m_{r,t}\end{pmatrix}=
              \begin{pmatrix}\widetilde{m}_1\\\widetilde{m}_2\\\vdots\\\widetilde{m}_r\end{pmatrix}$$
as characteristic invariants (hence, fulfilling properties (1)-(4) in Theorem~\ref{thm_characterization_cs}) and matrix
$$\mathcal{L}=\begin{pmatrix}l_{1,1}&l_{1,2}&\dots&l_{1,t}\\l_{2,1}&l_{2,2}&\dots&l_{2,t}\\
              \vdots&\vdots&\ddots&\vdots\\l_{s,1}&l_{s,2}&\dots&l_{s,t}\end{pmatrix}=
              \begin{pmatrix}\widetilde{l}_1\\\widetilde{l}_2\\\vdots\\\widetilde{l}_s\end{pmatrix}$$
of shift-maximal losing vectors if and only if the complete game $(N, \mathcal W^{\ast})$ has vector $\widetilde{n}=(n_1,\dots,n_t)\in\mathbb{N}_{>0}^t$ and  matrix
              $$\mathcal{M^{\ast}}=\begin{pmatrix}n_1-l_{s,1}&n_2-l_{s,2}&\dots&n_t-l_{s,t}\\n_1-l_{s-1,1}&
              n_2-l_{s-1,2}&\dots&n_t-l_{s-1,t}\\
              \vdots&\vdots&\ddots&\vdots\\n_1-l_{1,1}&n_2-l_{1,2}&\dots&n_t-l_{1,t}\end{pmatrix}=
              \begin{pmatrix}\widetilde{m_1^{\ast}}\\\widetilde{m_2^{\ast}}\\
              \vdots\\\widetilde{m_s^{\ast}}\end{pmatrix}$$
as characteristic invariants (hence, fulfilling properties (1)-(4) in Theorem~\ref{thm_characterization_cs}-(a)) and matrix
$$\mathcal{L}^{\ast}=\begin{pmatrix}n_1-m_{r,1}&n_2-m_{r,2}&\dots&n_t-m_{r,t}\\
n_1-m_{r-1,1}&n_2-m_{r-1,2}&\dots&n_t-m_{r-1,t}\\
              \vdots&\vdots&\ddots&\vdots\\n_1-m_{1,1}&n_2-m_{1,2}&\dots&n_t-m_{1,t}\end{pmatrix}=
              \begin{pmatrix}\widetilde{l_1^{\ast}}\\\widetilde{l_2^{\ast}}\\\vdots\\\widetilde{l_r^{\ast}}
              \end{pmatrix}$$
of shift-maximal losing vectors.
\end{itemize}
\end{lemma}

Let $\overline{cg}(n,t,r)$ be the number of complete games with $n$ players, with $t$ equivalence classes and $r$ shift-maximal losing vectors and similar notation for: $\overline{cg}(n,\ast,r)$, $\overline{wg}(n,t,r)$ and $\overline{wg}(n,\ast,r)$. Lemma~\ref{l:duals} allows us to deduce the next corollary and describes how  the bijection for characteristic invariants works.

\begin{corollary} \label{c:formulas dual}
${cg}(n,t,r) = \overline{cg}(n,t,r)$ and  ${wg}(n,t,r) = \overline{wg}(n,t,r)$.
\end{corollary}

The application of this Corollary and the results on enumerations in Sections 2 and 3 to games with one shift-maximal losing vector gives the enumerations for complete games and for weighted games respectively.
An analogous version of Theorem~\ref{t:cg} is given by
\begin{corollary}\label{c:cg}
\begin{enumerate}
\item $\overline{cg}(n,\ast,1)=2^n-1$,
\item $$\overline{cg}(n,t,1) = \left\{
                     \begin{array}{ll}
                       {n}, & \hbox{if} \ t=1 \\ \\
                       {n+1 \choose 2t-1}, & \hbox{if} \ 2 \leq t \leq \frac n2 + 1 \\ \\
                       {0}, & \hbox{otherwise}
                     \end{array}
                   \right.
$$
\end{enumerate}
\end{corollary}
and an analogous version of Corollary~\ref{c:compactwg(n,1)} is given by
\begin{corollary} \label{c:compactwgdual(n,1)}
$$ \overline{wg}(n,\ast,1) = \left\{
                        \begin{array}{ll}
                          2^n-1, & \hbox{if} \; n \leq 5 \\
                          \dfrac{n^4-6n^3+23n^2-18n+12}{12}, & \hbox{if} \; n \geq 6 \\
                        \end{array}
                      \right.
$$
\end{corollary}

Due to Proposition~\ref{p:ds79} and Lemma~\ref{l:duals} it follows that the results in Section 4 for the Shapley-Shubik index extend to complete games with one-shift maximal losing vector, since for any given $(n_1, n_2)$ and $\mathcal M = (a \ b)$, the dual game is given by the same vector and matrix $\mathcal L^{\ast} = (n_1-a \ n_2-b)$ which corresponds to matrix $\mathcal M^{\ast}$ where:
\begin{itemize}
\item[$(i)$]
$\mathcal M^{\ast} = (n_1-a+1 \ \ \  0)$ \; if $b=0$, otherwise:
\item[$(ii)$] $\mathcal M^{\ast} = \left( \begin{array}{cc} n_1-a+1 \ & \  0 \\ n_1-a-b+1 \ &\  n_2 \\ \end{array} \right) $ \; if $a+b-1 \leq n_1$,
\item[$(iii)$] $\mathcal M^{\ast} = \left( \begin{array}{cc} n_1-a+1 \ &\  0 \\ 0 \ &\  n-a-b+1 \\ \end{array} \right) $ \; if  $a+b-1 > n_1$.
\end{itemize}
Thus Proposition~\ref{p:SSpower} has a simple analogue, and therefore Corollary~\ref{cor_pi_orderSS} has a simple analogue too. The conjecture stated for the relative Banzhaf index derived for the computation up to $100$ players is also open in this dual context.

\section{Future work}

Any progress concerning enumeration of games, like $cg(n,t,r)$ or $wg(n,t,r)$, will be a significant advance in both directions: either providing new formulas or providing tighter bounds. In~\cite{KuTa12} it is proved that $cg(n,t,r)$ is a quasi-polynomial in $n$, if $t$ and $r$ are given, and can therefore be automatically computed (without escaping of the problem of capacity limitations). We wonder whether these automatic computations can be performed for the number $wg(n,t,r)$ of weighted games.

We also encourage research in other classes of simple games. For instance, roughly weighted games (considered in~\cite{TaZw99} and extensively studied in~\cite{GoSl11}) which are complete, i.e., roughly weighted complete games with $n$ players and $t$ types of equivalent players to determine $rcg(n,t,r)$.

As a future research concerning section 4 (and 5) it would be nice to find the proofs of our conjectures~\ref{c:BzR1} and~\ref{c:BzR2}. It seems likely that these conjectures are true since they hold for games with less than 100~players. It would also be of interest to know whether the Shapley-Shubik power index is the unique semiindex (semivalue restricted to simple games) that satisfies Corollary~\ref{cor_pi_orderSS}.

\section*{Acknowledgments}
The authors are grateful to the two referees of this paper for their
interesting comments and also for their exhaustive reports that contributed
to improve the original submitted version.

\providecommand{\bysame}{\leavevmode\hbox to3em{\hrulefill}\thinspace}
\providecommand{\MR}{\relax\ifhmode\unskip\space\fi MR }
\providecommand{\MRhref}[2]{%
  \href{http://www.ams.org/mathscinet-getitem?mr=#1}{#2}
}
\providecommand{\href}[2]{#2}

\end{document}